# A Geometric Identity for Pappus' Theorem

Michael Hawrylycz
Group C-3 Los Alamos National Laboratory
Los Alamos, New Mexico 87544

September 16, 1994

**Abstract**

An expression in the exterior algebra of a Peano space yielding Pappus' Theorem was originally given by Doubilet, Rota, and Stein[3]. Motivated by an identity of Rota, we give an identity in a Grassmann-Cayley algebra of step 3, involving joins and meets alone, which expresses the Theorem of Pappus.

While the bracket algebra itself has been shown to be a highly effective tool for projective space computations [1],[3],[7], Grassmann-Cayley algebra identities, those involving join and meet alone, have the advantage that the join and meet of extensors are again extensors and can be easily interpreted geometrically. Sturmfels and Whitely [5] have shown that any homogeneous bracket polynomial can be factored as a Grassmann-Cayley algebra expression, after multiplication by a suitable product of brackets. Their result, however, does not afford a practical means of generating identities in a Grassmann-Cayley algebra as the intial factor of brackets may be unwieldy. Although few general techniques are known for generating Grassmann-Cayley algebra identities, an elegant identity representing Desargues' Theorem has been given [1]. We give a new identity for the Theorem of Pappus, based on a similar identity in a Peano space due to Rota [4].

**Theorem 1** *In a Grassmann-Cayley algebra of step 3, the following identity is valid:*

$$(bc' \wedge b'c) \vee (ca' \wedge c'a) \vee (ab' \wedge a'b) = (c'b \wedge b'c) \vee (ca' \wedge ab) \vee (ab' \wedge a'c') \qquad (0.1)$$

The proof of Theorem 1 is a consequence of the following Proposition:

**Proposition 2** *If $a, b, c, a', b', c'$ are any six points in the plane then the expression*

$$(bc' \wedge b'c) \vee (ca' \wedge c'a) \vee (ab' \wedge a'b)$$

*changes at most in sign under any permutation of the points.*

*Proof.* We may first expand each of the parenthesized terms to write this expression as

$$([bc'c]b' - [bc'b']c) \vee ([ca'a]c' - [ca'c']a) \vee ([ab'b]a' - [ab'a']b)$$

Expanding further;

$$+[bc'c][ca'a][ab'b][b'c'a'] - [bc'c][ca'a][ab'a'][b'c'b]$$

$$-[bc'c][ca'c'][ab'b][b'aa'] + [bc'c][ca'c'][ab'a'][b'ab]$$

$$-[bc'b'][ca'a][ab'b][cc'a'] + [bc'b'][ca'a][ab'a'][cc'b]$$

$$+[bc'b'][ca'c'][ab'b][caa'] - [bc'b'][ca'c'][ab'a'][cab]$$

of which the only surviving terms are

$$[bc'c][ca'a][ab'b][b'c'a'] - [bc'b'][ca'c'][ab'a'][cab] \qquad (0.2)$$

Expression (0.2) changes sign alone under interchanging any of the points in the pairs $bc', b'c, ca', c'a, ab', a'b$. If the points $a, b$ are interchanged then the sum of (0.2) and the expression obtained may be written as

$$([aa'b'][bb'c'] - [ba'b'][ab'c'])[cc'a'][abc] + ([cac'][bca'] - [bcc'][caa'])[abb'][a'b'c'] \qquad (0.3)$$

By expanding the meet $aa'b \wedge b'$ in two different ways we obtain the identity $[aa'b']b + [a'bb']a + [bab']a' = [aa'b]b'$ which may be meeted with $b'c'$ to obtain

$$[aa'b'][bb'c'] + [a'bb'][ab'c'] + [bab'][a'b'c'] = 0 \qquad (0.4)$$

Substituting relation (0.4) into the sum (0.3), the first term in parentheses reduces to $[abb'][a'b'c']$. By similarly splitting the extensor $abc \wedge c'$ and meeting with $a'c$ we obtain

$$[cac'][bca'] + [cbc'][caa'] + [cc'a'][abc] = 0 \qquad (0.5)$$

which reduces the second parenthesized term to $-[cc'a'][abc]$ and thus the sum vanishes. By symmetry the expression changes sign under any transposition of unprimed or primed letters, and the proposition follows. □.

From Proposition 2 we may directly interpret the identity of Theorem 1.

**Theorem 3 (Pappus)** *If $a, b, c$ are collinear and $a', b', c'$ are collinear and all distinct, then the intersections $ab' \cap a'b$, $bc' \cap b'c$ and $ca' \cap c'a$ are collinear.*

*Proof.* The identity (0.1) follows from Proposition 2 by interchanging vectors $b$ and $c'$. To obtain Pappus theorem, assume $a, b, c$ and $a', b', c'$ are two sets of three collinear points in the plane. Now $ca' \wedge ab = c[a'ab]$ and $ab' \wedge a'c' = -b'[aa'c']$ and so the right hand side of (0.1) reduces to $(c'b \wedge b'c) \vee c \vee b' = 0$. Hence the left hand side vanishs as well. But the left hand side of (0.1) vanishs most generally when the points $bc' \wedge b'c, ca' \wedge c'a$, and $ab' \wedge a'b$ are collinear. $\square$

**Acknowledgement.** I am indebted to Gian-Carlo Rota for suggesting the problem and for several helpful discussions.